\documentclass[12pt]{article}
\usepackage{latexsym,amsmath,amssymb}
\usepackage{amsfonts}

\numberwithin{equation}{section} \numberwithin{table}{section}

\newcommand{\pd}[2]{\frac{\partial #1}{\partial #2}}

\newcommand{\vep}{\mbox{$\varepsilon$}}
\newcommand{\pa}{\mbox{$\partial$}}

\newcommand{\ode}{ordinary differential equation}
\newcommand{\pde}{partial differential equation}

\newcommand{\DE}{differential equation}

\newcommand{\PDE}{partial differential equation}

\newcommand{\LSA}{Lie symmetry analysis}
\def\d {\mbox {\rm d}}

\newcommand{\ds}{\displaystyle}

\begin{document}

\begin{center}
{\large \bf Conservation laws for a mathematical model of HIV transmission}\\[3 mm]

{ Winter Sinkala , Andrew Otieno\footnote{Corresponding author. Tel.: +27 47 502 2271; Fax: +27 47 502 2725.\\E-mail addresses: aotieno@wsu.ac.za , wsinkala@wsu.ac.za}}\\[3 mm]
Department of Mathematical Sciences and Computing,
Walter Sisulu University,  Private Bag X1, Mthatha 5117, South Africa.
\end{center}

\vspace{-5mm}

\begin{abstract}
\noindent A theorem due to Nail H. Ibragimov (2007) provides  a connection between symmetries and conservation
laws for arbitrary differential equations. The theorem is valid for
any system of differential equations provided that the number of
equations is equal to the number of dependent variables. In this
paper we use the theorem to determine conservation laws for a
nonlinear system of \ode s that represents a mathematical model for
HIV transmission.
\end{abstract}

\smallskip

\section{Introduction}
Ibragimov \cite{Ibragimov2007} provides a general theorem on conservation laws by means of which conservation laws can be constructed for an arbitrary system of \DE s admitting Lie symmetries, provided the number  of equations is equal to the number of dependent variables. Ibragimov's theorem is an extension of the classical theorem of Noether  \cite{Noether1918} in that the latter theorem does not require existence of a Lagrangian. In fact, unlike Noether's theorem, Ibragimov's theorem allows one to associate a conservation law to every Lie point symmetry admitted by a given arbitrary system of \DE s. Applications that one comes across involving the use of admitted Lie point symmetries to construct conservation laws via either the classical Noether's theorem or Ibragimov's theorem are often limited to scalar \PDE s \cite{Mahomed,Khalique,Wei&Zhang15,Freire2014,Govinder,Anco,Bluman,karaM,Yasara10}. The authors were unable to find any applications involving systems of first-order \ode s, for instance.
In this connection the application reported in this paper adds to the repertoire of nontrivial applications of Ibragimov's theorem.

We consider a nonlinear system of three first-order ordinary differential equations that arise from a model formulated by Anderson~\cite{Anderson} which describes the transmission of HIV/AIDS in male homosexual/bisexual cohorts. A particular case of the model translates into a coupled nonlinear system of first-order ordinary differential equations \cite{Torrisi}:
\begin{equation}\label{gptra4abc}
\begin{array}{l}
 \ds F_1(t,u,u_{(1)})  =  u^1_t +\frac{\beta cu^1u^2}{u^1+u^2+u^3}+\mu
 u^1= 0\\
 \ds F_2(t,u,u_{(1)})   = u^2_t -\frac{\beta
 cu^1u^2}{u^1+u^2+u^3}+(\nu+\mu)u^2 = 0\\
 \ds F_3(t,u,u_{(1)})   = u^3_t - \nu u^2 + (\mu+\beta c) u^3 = 0,
\end{array}
\end{equation}
where $u_{(1)}$ denotes the first-order partial derivatives $u^1_t$, $u^2_t$ and $u^3_t$.
This system is typical of models that are formulated to mimic dynamics in epidemiology. Such systems are often highly nonlinear and difficult to analyze.
In their seminal work, Torrisi and Nucci \cite{Torrisi} perform Lie symmetry analysis on (\ref{gptra4abc}). They determine a solvable Lie algebra admitted by the system and exploit this to find a solution (albeit in quadrature form) of the system. In this article  we extend the Lie symmetry analysis of (\ref{gptra4abc}) by generating conservation laws for the system via Ibragimov's new conservation theorem  \cite{Ibragimov2007}.

Conservation laws of physical systems are fundamental to our understanding of the system being studied. Apart from having a direct physical interpretation, conservation laws may be essential in studying the integrability of the system. For example, in the numerical integration of \pde s conservation laws help to control numerical errors in that they describe the properties of the system that do not change. On the whole conservation laws play an important role in the analysis of basic properties of the solutions. Therefore, the construction of conservation laws is one of the most important applications of symmetries to physical systems \cite{olver,BlumanKumei}.

The rest of this paper is organised as follows. We present elements of Lie symmetry analysis of \DE s in Section~\ref{Sec2}. An overview of the theorems of Noether and Ibragimov for constructing conservation laws via admitted symmetries is provided in Section~\ref{Sec3}. Conservation Laws of Anderson's HIV model are constructed in Section~\ref{Sec4}. In Section~\ref{Sec4} we discuss the results and give concluding remarks.

\section{Preliminaries}\label{Sec2}
Let
us consider an $r$th-order ($r\ge 1$) system of \DE s with  $m$ dependent
variables $u=(u^1,u^2,\ldots,u^m)$ and $n$ independent variable $x=(x^1,x^2,\ldots,x^n)$, $u=u(x)$,
\begin{equation}\label{3.1}
    F_\alpha(x,u,\ldots,u_{(r)}) = 0, \quad \alpha = 1,\ldots,m,
\end{equation}
where $u_{(k)}$ denotes the collection $\left\{u_k^\alpha\right\}$ of $k$th-order derivatives, $k\ge 1$.
Suppose that (\ref{3.1}) admits a one parameter Lie group of point transformations
\begin{equation}\label{3.2}
    \begin{array}{c}
       \widetilde{x}^i =x^i +\vep \,\xi^i(x,u)+O(\vep^{2}),\\
       \widetilde{u}^\alpha =u^\alpha +\vep\,\eta^\alpha(x,u)+O(\vep^{2}), \\
\end{array}
\end{equation}
where $\vep$ is a real parameter;  $\xi^i$ and $\eta^\alpha$ are given smooth functions.
Invariance of (\ref{3.1}) under (\ref{3.2}) is conveniently
expressed in terms of the infinitesimal generator of (\ref{3.2}) is
\begin{equation}\label{3.3}
X =
\xi^i(x,u)\dfrac{\partial}{\partial{x^i}}+\eta^\alpha(x,u)\dfrac{\partial}{\partial{u^\alpha}},
\end{equation}
where the usual convention of summation over repeated indices is adopted \cite{BlumanKumei}. In fact this convention is adopted in all subsequent expressions.
We say that (\ref{3.2}) is a symmetry of (\ref{3.1}) if and only if for
all $\alpha = 1,\ldots,m$
\begin{equation}\label{3.4}
    X^{(r)}[F_\alpha(x,u,\ldots,u_{(r)})]\big\vert_{(\ref{3.1})} =
    0,
\end{equation}
where $X^{(r)}$ is the $r$th extension of (\ref{3.3}) defined by
\begin{eqnarray}
 X^{(r)} &=& \xi^i(x,u)\pa_{x^i} + \eta^\alpha(x,u)\pa_{u^\alpha} +
 \zeta_i^{\alpha}(x,u,u_{(1)})\pa_{u_i^\alpha} + \zeta_{i_1i_2}^{\alpha}(x,u,u_{(1)},u_{(2)})\pa_{u_{i_1i_2}^\alpha}\nonumber\\
 && +\, \cdots +
 \zeta_{i_1i_2 \ldots i_r}^{\alpha}(x,u,u_{(1)},\ldots,u_{(r)}) \pa_{u_{i_1i_2 \ldots i_r}^\alpha},
\end{eqnarray}
with the explicit formulas for the extended infinitesimal coefficients given
recursively by
\begin{eqnarray}
 \zeta_i^{\alpha} &=& D_i\left(\eta^\alpha\right) - u_j^\alpha D_i\left(\xi^j\right),\quad i=1,2,\ldots,n,\\
\zeta_{i_1i_2 \ldots i_{s}}^{\alpha} &=& D_{i_s}\left(\zeta_{i_1i_2 \ldots
i_{s-1}}^{\alpha}\right)- u_{j i_1 i_2 \ldots i_{s-1}}^{\alpha} D_{i_s}\left(\xi^j\right), \quad s > 1,
\end{eqnarray}
where $D_i$
is the total derivative operator with respect
to $x^i$ defined by
\begin{equation}\label{2.8}
D_i = \frac{\pa}{\pa x^i} + u_i^\alpha \frac{\pa}{\pa u^{\alpha}} + u_{ij}^\alpha \frac{\pa}{\pa u_j^{\alpha}} + \cdots+ u_{i i_1 i_2\ldots i_n}^\alpha \frac{\pa}{\pa u_{i_1 i_2\ldots i_n}^\alpha} + \cdots.
\end{equation}
Note that in terms of the total derivative operator the derivatives of $u^\alpha$ with respect to $x_i$ are
 \[u_i^\alpha = D_i\left(u^\alpha\right), \quad u_{ij}^\alpha = D_i\left(u_j^\alpha\right) - D_i D_j\left(u^\alpha\right), \quad \ldots\]
A conserved vector of (\ref{3.1}) is an $n$-tuple
\begin{equation}\label{5.2}
    C =\left(C^1(x,u,\ldots,u_{(r-1)}),\ldots, C^n(x,u,\ldots,u_{(r-1)})\right)
\end{equation}
such that
\begin{equation}\label{2.10}
    D_i\left(C^i\right) = 0
\end{equation}
on the solution space of (\ref{3.1}). The expression (\ref{2.10}) is a conservation law of (\ref{3.1}).

\section{The connection between conservation laws and admitted symmetries}\label{Sec3}
A fundamental relationship between symmetries and conservation laws is provided by  Noether's theorem \cite{Noether1918}, which states that for Euler-Lagrange systems of differential equations, to each \emph{Noether symmetry} associated with the Lagrangian there corresponds a conservation law which can be determined explicitly by a formula \cite{BlumanKumei}. Noether's theorem therefore reduces the search for conservation laws to a search for Noether symmetries. However, the dependence upon knowledge of a suitable Lagrangian to exploit Noether's theorem diminishes the applicability of the theorem significantly. Ibragimov's theorem \cite{Ibragimov2007} extends the application of Noether's theorem by providing for the association of a conservation law to every symmetry of a system of differential equations, albeit with the proviso that the number of equations in the system equals the number of dependent variables and that the given system be considered together with the associated adjoint system. The rest of this section introduces the essential elements of Noether's theorem and Ibragimov's theorem.

 \subsection{Noether's theorem~\cite{Noether1918}}
Consider a system of \DE s identical with Euler-Lagrange equations
\begin{equation}\label{4.1}
    \frac{\delta L}{\delta u^\alpha}\equiv \pd{L}{u^\alpha} - D_i\left(\pd{L}{u_i^\alpha}\right) = 0,\quad \alpha=1,\ldots,m,
\end{equation}
arising from the variational integral
\begin{equation}\label{4.2}
    \mbox{$\int L(x,u,\ldots,u_{(r)}) \d x$}
\end{equation}
taken over an arbitrary $n$-dimensional domain in the space of the independent variables $x=(x^1,\ldots,x^n)$. The Lagrangian $L(x,u,\ldots,u_{(r)})$ involves $x$ and the dependent variables $u=(u^1,\ldots,u^n)$, $u = u(x)$, together with their derivatives $u_{(i)}$.

Let the system of \DE s (\ref{3.1}) admit a continuous group $G$ with a generator (\ref{3.3}).
Noether's theorem states that if the variational integral (\ref{4.2}) is invariant under the group $G$, then the vector field $C =\left(C^1,\ldots, C^n\right)$ defined by
\begin{eqnarray}\label{ik4}
  C^i &=& L \xi^i + \left(\eta^\alpha - \xi^i u_j^\alpha\right)\left[\pd{L}{u_i^\alpha} - D_j\left(\pd{L}{u_{ij}^\alpha}\right) + D_jD_k\left(\pd{L}{u_{ijk}^\alpha}\right) - \cdots\right]\nonumber \\[0.6ex]
    &&+\,\,  D_j\left(\eta^\alpha - \xi^i u_j^\alpha\right)\left[\pd{L}{u_{ij}^\alpha} - D_k\left(\pd{L}{u_{ijk}^\alpha}\right) +\cdots\right]\nonumber \\[0.6ex]
  && +\,\, D_jD_k \left(\eta^\alpha - \xi^i u_j^\alpha\right)\left[\pd{L}{u_{ijk}^\alpha} - \cdots\right] + \cdots \label{3.16}
\end{eqnarray}
 provides a conservation law for the Euler-Lagrange equations (\ref{4.1}). Noether's theorem states that if the Invariance of the variational integral (\ref{4.2}) under the group $G$ is established via the infinitesimal test for invariance,
\begin{equation}\label{4.5}
  X\left(L\right) + L D_i\left(\xi^i\right) = 0,
\end{equation}
where the appropriate prolongation of $X$ is understood.

\subsection{Extension of Noether's theorem: Conservation Laws via Ibragimov's theorem~\cite{Ibragimov2007}}

Consider a system of $r$th-order differential equations defined in (\ref{3.1}).
We introduce the differential
functions
\begin{equation}\label{}
  F_\alpha^*(x,u,v,\ldots,u_{(r)},v_{(r)}) = \frac{\delta\left(v^\beta F_\beta\right)}{\delta u^\alpha}, \quad
  \alpha=1,\ldots,m
\end{equation}
where $v = (v^1,\ldots, v^m)$ are new dependent variables, $v = v(x)$, and
\begin{equation}\label{3.8}
  \frac{\delta}{\delta u^\alpha} = \pd{}{u^\alpha} + \sum (-1)^r
  D_{i_1},\ldots,D_{i_r}\pd{}{u_{i_1,\ldots,i_r}^\alpha},\quad
  \alpha=1,\ldots,m.
\end{equation}
The system of adjoint equations to (\ref{3.1}) is defined by
\begin{equation}\label{5.6}
  F_\alpha^*(x,u,v,\ldots,u_{(r)},v_{(r)}) = 0, \quad
  \alpha=1,\ldots,m.
\end{equation}
We now have that the simultaneous system consisting of the $r$th-order differential equations (\ref{3.1}) considered together with its
adjoint equation (\ref{5.6}) has a Lagrangian defined by
\begin{equation}\label{5.7}
  L = v^\beta F_\beta(x,u,\ldots,u_{(r)}).
\end{equation}
Furthermore, the adjoint system (\ref{5.6}) inherits the symmetries of the system (\ref{3.1}) in the sense that if the system (\ref{3.1}) admits a point transformation
group with a generator
\begin{equation}\label{5.8}
  X = \xi^i(x,u)\pd{}{x^i} + \eta^\alpha(x,u)\pd{}{u^\alpha},
\end{equation}
then the adjoint system (\ref{5.6}) admits the operator (\ref{5.8}) extended to the variables $v^\alpha$ by the
formula
\begin{equation}\label{5.9}
  Y = \xi^i\pd{}{x^i} + \eta^\alpha\pd{}{u^\alpha} + \eta_*^\alpha\pd{}{v^\alpha},
\end{equation}
with coefficients $\eta_*^\alpha = \eta_*^\alpha(x,u,v,\ldots)$ chosen in such a way that $Y$ satisfies the infinitesimal test for invariance of the variational integral associated with (\ref{5.7}), i.e,
\begin{equation}\label{}
  Y\left(L\right) + L D_i\left(\xi^i\right) = 0,
\end{equation}
where the generator $Y$ is prolonged appropriately to the $r$th derivatives $u_{(r)}$ and $v_{(r)}$. It turns out that
\begin{equation}\label{3.14}
  \eta_*^\alpha = - \left[\lambda_\beta^\alpha v^\beta + v^\alpha D_i\left(\xi^i\right)\right], \quad
  \alpha=1,\ldots,m.
\end{equation}
with $\lambda_\beta^\alpha$ defined by the invariance condition (\ref{3.4}) in the form
\begin{equation}\label{3.15}
  X\left(F_\alpha\right) = \lambda_\beta^\alpha F_\beta, \quad
  \alpha=1,\ldots,m,
\end{equation}
where the prolongation of $X$ to all derivatives involved in Lagrangian the system (\ref{3.1}) is understood.

Noether's theorem is now employed to furnish the conserved vector $C^i =\left(C^1,\ldots, C^n\right)$, with $C_i$ defined by
\begin{eqnarray}
  C^i &=& L \xi^i + \left(\eta^\alpha - \xi^i u_j^\alpha\right)\left[\pd{L}{u_i^\alpha} - D_j\left(\pd{L}{u_{ij}^\alpha}\right) + D_jD_k\left(\pd{L}{u_{ijk}^\alpha}\right) - \cdots\right]\nonumber \\[0.6ex]
    &&+\,\,  D_j\left(\eta^\alpha - \xi^i u_j^\alpha\right)\left[\pd{L}{u_{ij}^\alpha} - D_k\left(\pd{L}{u_{ijk}^\alpha}\right) +\cdots\right]\nonumber \\[0.6ex]
  && +\,\, D_jD_k \left(\eta^\alpha - \xi^i u_j^\alpha\right)\left[\pd{L}{u_{ijk}^\alpha} - \cdots\right] + \cdots \label{3.16}
\end{eqnarray}
where $u^\alpha = \left(u^1,\ldots,u^m,v^1,\ldots,v^m \right).$
\section{Conservation Laws of Anderson's HIV model}\label{Sec4}
For the problem at hand, (\ref{gptra4abc}), we have a system of three dependent
variables $u=(u^1,u^2,u^3)$ and one independent variable $t$, where $u=u(t)$.
This system is considered together with the corresponding adjoint system of equations which is constructed as outlined in Section~\ref{Sec3}.  Let $v=(v^1,v^2,v^3)$ be the new dependent variables, $v=v(t)$. According to (\ref{5.6}) the adjoint system is given by
\begin{equation}\label{adjt2}
  F_\alpha^*(t,u,v,u_{(1)},v_{(1)}) = \frac{\delta\left(v^\beta F_\beta\right)}{\delta u^\alpha}=0, \quad
  \alpha=1,2,3,
\end{equation}
where $\delta/\delta u^\alpha$ is the Euler-Lagrange operator defined by (\ref{3.8}), which for the problem being considred reduces to
\begin{equation}\label{}
  \frac{\delta}{\delta u^\alpha} = \pd{}{u^{\alpha}}  -
  D_t\pd{}{u_{t}^{\alpha}},\quad
  \alpha=1,2,3,
\end{equation}
where $D_t$
is the total derivative operator with respect to $t$ defined by
\[D_t = \frac{\pa}{\pa t} + u^1_t\frac{\pa}{\pa u^1} + u^2_t\frac{\pa}{\pa u^2} + u^3_t\frac{\pa}{\pa u^3} + \cdots.\]

Thus (\ref{adjt2}) translates into the following adjoint system of equations:
\begin{eqnarray}\label{6.5}
 \frac{dv^1}{dt}&=& v^1\left(\frac{u^2\,\beta \,\delta }{u^1 + u^2 + u^3} - \frac{u^1\,u^2\,\beta \,\delta }{{\left( u^1 + u^2 + u^3 \right) }^2} + \mu  \right) -\frac{v^2\,u^2\left( u^2 + u^3 \right) \,\beta \,\delta }{{\left( u^1 + u^2 + u^3 \right) }^2},\nonumber\\
 \frac{dv^2}{dt}&=& v^2\left( \mu  + \nu - \frac{u^1\left( u^1 + u^3 \right) \,\beta \,\delta }{{\left( u^1 + u^2 + u^3 \right) }^2}\right) + \frac{v^1\,u^1\left( u^1 + u^3 \right) \,\beta \,\delta }{{\left( u^1 + u^2 + u^3 \right) }^2} - v^3\,\nu,\\
 \frac{dv^3}{dt}&=&v^3\,\alpha  - \frac{v^1\,u^1\,u^2\,\beta \,\delta }{{\left( u^1 + u^2 + u^3 \right) }^2} + \frac{v^2\,u^1\,u^2\,\beta \,\delta }{{\left( u^1 + u^2 + u^3 \right) }^2}\nonumber.
\end{eqnarray}
According to the infinitesimal condition for invariance (\ref{3.4}) the system (\ref{gptra4abc}) admits a symmetry group with the infinitesimal generator
\[X=\xi(t,u^1,u^2,u^3)\partial_{t}+\eta^1(t,u^1,u^2,u^3)\partial_{u^1}+\eta^2(t,u^1,u^2,u^3)
\partial_{u^2}+\eta^3(t,u^1,u^2,u^3)\partial_{u^3}\]
if and only if
\begin{equation}\label{gptra22abc}
\left. X^{(1)}F_1\right\vert_{(\ref{gptra4abc})}=0,\quad
\left. X^{(1)}F_2\right\vert_{(\ref{gptra4abc})}=0,\quad
\left. X^{(1)}F_3\right\vert_{(\ref{gptra4abc})}=0,
\end{equation}
where
\begin{equation}\label{ext1}
X^{(1)}\!=\!\xi\partial_{t}+\eta^1\partial_{u^1}+\eta^2
\partial_{u^2}+\eta^3\partial_{u^3}+\zeta_t^1\partial_{u_t^1}+\zeta_t^2\partial_{u_t^2}+\zeta_t^3\partial_{u_t^3}
\end{equation}
with
\begin{eqnarray}
  \zeta_t^1 &=& D_t\left(\eta^1\right) - u^1_t D_t\, \xi\label{eqns 114a}\\
  \zeta_t^2 &=& D_t\left(\eta^2\right) - u^2_t D_t\, \xi\label{eqns 114b}\\
  \zeta_t^3 &=& D_t\left(\eta^3\right) - u^3_t D_t\, \xi.\label{eqns 114c}
\end{eqnarray}
After making an ansatz on the form of the operator $X$ (assuming
that the functions $\xi$, $\eta^{1}$ and $\eta^{2}$ are polynomials
of second
 degree on $u^{1}$, $u^{2}$ and $u^{3}$) the solution of the equations (\ref{gptra22abc}), after lengthy analysis, leads to a three-dimensional
 Lie symmetry algebra $L_3 = \left<X_1, X_2, X_3\right>$ admitted by (\ref{gptra4abc}) with the following basis operators \cite{Torrisi}:
%
%
\begin{eqnarray}\label{6.4}
X_1 &=& \partial_t,\nonumber\\[1.3ex]
X_2 &=& u^1\,\partial_{u^1}+u^2\,\partial_{u^2}+u^3\,\partial_{u^3},\\[1.3ex]
X_3 &=&
e^{-(\mu+\nu)\,t}\partial_{u^2}+\frac{u^1+u^3}{u^2}\,\,e^{-(\mu+\nu)\,t}\,
\partial_{u^3}.\nonumber
\end{eqnarray}

We shall find conservation laws for the simultaneous system (\ref{gptra4abc}) and (\ref{6.5}) using each of the symmetries in (\ref{6.4}), extended suitably to the adjoint variables $v^{1}$, $v^{2}$ and $v^{3}$. According to (\ref{5.7}) the system (\ref{gptra4abc}) and (\ref{6.5}) considered together has the Lagrangian
\begin{eqnarray}\label{4.15}
    L &=& v^\beta F_\beta =  v^1\left[u_t^1 - \frac{\beta cu^1u^2}{u^1+u^2+u^3}+\mu
 u^1\right] - v^2\left[u_t^2 - \frac{\beta
 c u^1u^2}{u^1+u^2+u^3}+(\nu+\mu)u^2\right]\nonumber\\ &&\quad\quad\quad +\, v^3\left[u_t^3 - \nu u^2 + (\mu+\beta c) u^3\right].
\end{eqnarray}
Each of the symmetries in (\ref{6.4}) has the form
\begin{equation}\label{}
X=\xi\, \partial_{t}+\eta^1 \partial_{u^1}+\eta^2 \partial_{u^2}+\eta^3 \partial_{u^3}
\end{equation}
and needs to be extended appropriately to the operator
\begin{equation}\label{4.17}
Y = X + \eta_*^1\,\partial_{v^1} + \eta_*^2\,\partial_{v^2} + \eta_*^3\,\partial_{v^3},
\end{equation}
to cater for the adjoint differential variables. It turns out that (\ref{4.17}) is admitted by the adjoint system (\ref{6.5})  provided the infinitesimal coefficients $\eta_*^\alpha$ in the operator are prescribed as follows \cite{Ibragimov2007}:
\begin{equation}\label{}
  \eta_*^\alpha = - \left[\lambda_\beta^\alpha v^\beta + v^\alpha D_t\left(\xi\right)\right], \quad \alpha=1,2,3,
\end{equation}
where  $\lambda_\beta^\alpha$ are defined by the invariance condition (\ref{3.15}). Taking $X_1$ from (\ref{6.4}) and extending it once to $X_1^{(1)} = X_1$\,, condition (\ref{3.15}) leads to a set of equations,
\begin{eqnarray}\label{}
 X_1^{(1)}F_1 &=& \lambda_1^1\,F_1 + \lambda_2^1,F_2 + \lambda_3^1\,F_3\\
 X_1^{(1)}F_2 &=& \lambda_1^2\,F_1 + \lambda_2^2,F_2 + \lambda_3^2\,F_3\\
 X_1^{(1)}F_3 &=& \lambda_1^3\,F_1 + \lambda_2^3,F_2 + \lambda_3^3\,F_3,
\end{eqnarray}
which must be solved for $\lambda_\beta^\alpha$. Clearly $X_1^{(1)}F_\alpha = 0$ for all $\alpha =  1,2,3$. It follows, therefore, that $\lambda_\beta^\alpha=0$ for all $\alpha$ and $\beta$. Furthermore, $\xi=1$ in the generator $X_1$, which leads us to the conclusion that in this case $\eta_*^\alpha = 0$ for all $\alpha =  1,2,3$. Therefore
\begin{eqnarray}
 Y_1 &=& X_1 = \partial_t.\label{4.17w}
\end{eqnarray}
For $X_2$ and $X_3$, the coefficient $\xi$ equals zero, which means that the second term in (\ref{3.14}) vanishes. Proceeding as we did for $X_1$ we determine $\lambda_\beta^\alpha$ in each of the two cases and obtain $Y_2$ and $Y_3$, via the desired extensions of $X_2$ and $X_3$ respectively:
\begin{eqnarray}\label{}
Y_2 &=& u^1\,\partial_{u^1}+u^2\,\partial_{u^2}+u^3\,\partial_{u^3} - v^1\,\partial_{v^1} - v^2\,\partial_{v^2} - v^3\,\partial_{v^3},\label{4.18w}\\[1.3ex]
Y_3 &=& e^{-(\mu+\nu)\,t}\left[\partial_{u^2}+\frac{u^1+u^3}{u^2}\,\partial_{u^3} - c\left(\frac{1}{u^2}\,\partial_{v^1}  - \frac{u^1 + u^3}{(u^2)^2}\,\partial_{v^2} +\frac{1}{u^2}\,\partial_{v^3}\right)\right].\label{4.19w}
\end{eqnarray}
By applying Noether's theorem to each of the generators $Y_i$ with the associated Lagragian (\ref{4.15})
we wish to find the corresponding conservation laws. Let us rename the dependent variables $v^1$, $v^2$ and $v^3$ of the adjoint system as $u^4$, $u^5$ and $u^6$, respectively, so that each of the generators (\ref{4.17w}), (\ref{4.18w}) and (\ref{4.19w}) is in
the form
\begin{equation}\label{4.21}
  Y = \xi(t,u)\pd{}{t} + \eta^\alpha(t,u)\pd{}{u^\alpha},
\end{equation}
where \[u = \left(u^1,u^2,u^3,v^1,v^2,v^3 \right) = \left(u^1,u^2,u^3,u^4,u^5,u^6 \right).\]
According to (\ref{2.10}) and (\ref{3.16}), the conservation law corresponding to (\ref{4.21}) is
\begin{equation}\label{}
  D_t\left\{\xi L + \left(\eta^\alpha - \xi u^\alpha\right)\left[\pd{L}{u^{\alpha}}  -
  D_t\pd{L}{u^{\alpha}}\right]\right\} = 0.
\end{equation}
We therefore obtain the following conservation laws:
\begin{eqnarray}\label{6.9}
D_t\left\{ \frac{u^1\,u^2\,\beta \,\delta\left( v^1 - v^2 \right)}{u^1 + u^2 + u^3} + v^1\,u^1\,\mu  + v^2\,u^2\left( \mu  + \nu  \right)  + v^3\left( u^3\,\alpha  - u^2\,\nu  \right)\right\} &=& 0\nonumber\\[1.3ex]
D_t\left\{ v^1\,u^1 + v^2\,u^2 + v^3\,u^3\right\} &=& 0\\[1.3ex]
D_t\left\{ e^{-(\mu+\nu)\,t}\left[v^2 + \frac{v^3\left( u^1 + u^3 \right) }{(u^1)^2}\right]\right\} &=& 0, \nonumber
\end{eqnarray}
corresponding to each of the extended symmetries $Y_1$, $Y_2$ and $Y_3$, respectively, where $v=\left(v^1,v^2,v^3\right)$ is any solution of the adjoint system of equations (\ref{6.5}).

\section{Concluding remarks}\label{Sec5}
  In this paper we have considered a nontrivial nonlinear system of first-order \ode s arising from a mathematical model formulated by Anderson \cite{Anderson} to describe the transmission of HIV/AIDS in male homosexual/bisexual cohorts. We have applied  Ibragimov's theorem and constructed conservation laws of the system considered together with the associated adjoint system. We have, however, not attempted to attach physical meaning to the adjoint system (\ref{6.5}) and/or the constructed conservation laws (\ref{6.9}). We defer this to possible future work on the model (\ref{gptra4abc}).
The application reported in this paper is an instructive application of Ibragimov's new theorem and may be used to construct conservation laws in other settings involving systems of first-order \ode s.

\section{Acknowledgment}

Support from the Directorate of Research Development of Walter Sisulu University is gratefully acknowledged.

\end{document}